\documentclass{article}

\usepackage{amssymb}

\def\ccc{\mathbb{C}}
\def\ord{\mathrm{ord}\,}
\def\I{\mathrm{I}}
\def\m{\mathrm{m}}
\def\rank{\mathrm{rank}}

\title{Semicontinuity of the \L{}ojasiewicz exponent}
\author{Arkadiusz P\l{}oski}
\date{January 2011}
\newtheorem{thm}{Theorem}[section]
\newtheorem{ex}[thm]{Example}
\newtheorem{lem}[thm]{Lemma}
\newtheorem{prop}[thm]{Proposition}
\newtheorem{rem}[thm]{Remark}

\begin{document}
\maketitle

\noindent{\bf Abstract.} We prove that the \L{}ojasiewicz exponent $l_0(f)$ of a~finite holomorphic germ $f:(\ccc^n,0)\to(\ccc^n,0)$ is lower semicontinuous in any multiplicity-constant deformation of $f$.\footnotetext{\noindent{}2010 Mathematics Subject Classification. Primary: 32S05; Secondary 14B05.}

\section{Introduction}

\sloppy Let $\ccc\{z\}$ denote the ring of convergent power series in $n$ variables $z=(z_1,\ldots,z_n)$. Any sequence of convergent power series $h=(h_1,\ldots,h_p)\in\ccc\{z\}^p$ without constant term defines the germ of a~holomorphic mapping $h:(\ccc^n,0)\to (\ccc^p,0)$. We put $\displaystyle\ord h=\inf_k\{\ord h_k\}$, where $\ord h_k$ denotes the order of vanishing of $h_k$ at $0$ (by convention $\ord 0=+\infty$). If $\displaystyle|\underbar{z}|=\max_{j=1}^n|\underbar{z}_j|$ for $\underbar{z}=(\underbar{z}_1,\ldots,\underbar{z}_n)\in\ccc^n$ then $\ord h$ for $h\neq 0$ is the largest $\alpha>0$ such that $|h(\underbar{z})|\leqslant c|\underbar{z}|^\alpha$ with a~constant $c>0$ for $\underbar{z}\in\ccc^n$ close to $0\in\ccc^n$.

Let $f=(f_1,\ldots,f_n)\in\ccc\{z\}^n$, $f(0)=0$, define a~finite holomorphic germ $f:(\ccc^n,0)\to(\ccc^n,0)$; i.e. such that $f$ has an isolated zero at the origin $0\in\ccc^n$ and let $\I(f)$ be the ideal of $\ccc\{z\}$ generated by $f_1,\ldots, f_n$. Then $\I(f)$ is of finite codimension in $\ccc\{z\}$ and the multiplicity $\m_0(f)$ of $f$ is equal by definition to $\dim_\ccc\,^{\ccc\{z\}}\!/\!_{\I(f)}$. There exist arbitrary small neighbourhoods $U$ and $V$ of $0\in\ccc^n$ such that the mapping $U\ni\underbar{z}\to f(\underbar{z})\in V$ is an $\m_0(f)$-sheeted branched covering, see \cite{L}, chapter 5, \S{}2.

Another important characteristic of a~finite germ $f:(\ccc^n,0)\to(\ccc^n,0)$ introduced and studied by M.~Lejeune-Jalabert and B.~Teissier in 1973--1974 seminar at the Ecole Polytechnique (in a~very general setting), see \cite{LJ-T}, is the \L{}ojasiewicz exponent $l_0(f)$ defined to be the smallest $\theta>0$ such that there exist a~neighbourhood $U$ of $0\in\ccc^n$ and a~constant $c>0$ such that
$$
|f(\underbar{z})|\geqslant c|\underbar{z}|^\theta\quad\mbox{for all }\underbar{z}\in U.
$$
The \L{}ojasiewicz exponent can be calculated by means of analytic arcs (see \cite{LJ-T}, \S{}5 and \cite{P2}, \S{}2) $\phi(s)=(\phi_1(s),\ldots,\phi_n(s))\in\ccc\{s\}^n$, $\phi(0)=0$, $\phi(s)\neq 0$ in $\ccc\{s\}^n$:
$$
l_0(f)=\sup_\phi\left\{\frac{\ord f\circ\phi}{\ord\phi}\right\}.
$$
The following lemma \cite{P1}, Corollary 1.4 will be useful for us.
\begin{lem}\label{lem1.1}
Let $f:(\ccc^n,0)\to(\ccc^n,0)$ be a~finite holomorphic germ. Then $l_0(f)\leqslant \m_0(f)$ with equality if and only if $\rank\left(\frac{\partial f_i}{\partial z_j}(0)\right)\geqslant n-1$.
\end{lem}
Now, let $h\in\ccc\{z\}$, $h(0)=0$, be a~convergent power series defining an isolated singularity at $0\in\ccc^n$ i.e. such that the gradient of $h$, $\nabla h=\left(\frac{\partial h}{\partial z_1},\ldots,\frac{\partial h}{\partial z_n}\right):(\ccc^n,0)\to(\ccc^n,0)$ is finite at $0\in\ccc^n$. Then $\mu_0:=\m_0(\nabla h)$ is the Milnor number of the singularity $h=0$. Teissier calculated in \cite{T1} $\mathcal{L}_0(h):=l_0(\nabla h)$ in terms of polar invariants of the singularity and proved that the \L{}ojasiewicz exponent $\mathcal{L}_0(h)$ is lower semicontinuous in any $\mu$-constant deformation of the singularity $h=0$. He showed also that if we don't assume $\mu\,$=\,constant that $\mathcal{L}_0(h)$ is neither upper or lower semicontinuous, see \cite{T2}. The ``jump phenomena'' of the \L{}ojasiewicz exponent was rediscovered by some authors, see \cite{Mc-N}. The aim of this note is to prove that the \L{}ojasiewicz exponent is lower semicontinuous in any multiplicity-constant deformation of the finite holomorphic germ. The proof is based on the formula for the \L{}ojasiewicz exponent given by the author in \cite{P2} (see also Lemma \ref{lem3.3} in Section \ref{sect3}).

\section{Result}

Let $f=(f_1,\ldots,f_n)\in\ccc\{z\}^n$, $f(0)=0$, define a~finite holomorphic germ. A~sequence $F=(F_1,\ldots,F_n)\in\ccc\{t,z\}^n$ of convergent power series in $k+n$ variables $(t,z)=(t_1,\ldots,t_k,z_1,\ldots,z_n)$ is a~deformation of $f$ if $F(0,z)=f(z)$ in $\ccc\{z\}$ and $F(t,0)=0$ in $\ccc\{t\}$. Then the sequence $(t,F(t,z))\in\ccc\{t,z\}^{k+n}$ defines a~holomorphic germ $(\ccc^{k+n},0)\to(\ccc^{k+n},0)$ of multiplicity $\m_0(f)$. Indeed, it is easy to check that the algebras $^{\ccc\{z\}}\!/_{\I(f)}$ and $^{\ccc\{t,z\}}\!/_{\I(t,F)}$ are $\ccc$-isomorphic.

We put $F_{\underbar{t}}=F(\underbar{t},z)\in\ccc\{z\}^n$ for $\underbar{t}\in\ccc^k$ close to $0$. Then $F_{\underbar{t}}(0)=0$ and $\m_0(F_{\underbar{t}})\leqslant\m_0(F_0)=\m_0(f)$ for $\underbar{t}\in\ccc^k$ close to $0$, see \cite{T}, chapter 2, \S{}5. We say that $F$ is a~multiplicity-constant deformation of the germ $f:(\ccc^n,0)\to(\ccc^n,0)$ if $\m_0(F_{\underbar{t}})=\m_0(F_0)$ for $\underbar{t}$ close to $0$.

The main result of this note is 
\begin{thm}\label{thm2.1}
Let $F\in\ccc\{t,z\}^n$ be a~multiplicity-constant deformation of the germ $f:(\ccc^n,0)\to(\ccc^n,0)$. Then
$$
l_0(F_0)\leqslant l_0(F_{\underbar{t}})\quad\mbox{ for } \underbar{t}\in\ccc^k\mbox{ close to }0\in\ccc^k.
$$
Moreover, if $F$ is a~one-parameter deformation ($k=1$), then $l_0(F_{\underbar{t}})$ is constant for $\underbar{t}\neq 0$ close to $0\in\ccc$. 
\end{thm}
The proof of the theorem is given in Section \ref{sect4} of this note. The inequality stated above may be strict:
\begin{ex}\rm (see \cite{Mc-N}, \S{}5).\\
Let $F(t,z_1,z_2)=(tz_1+z_1^a+z_2^b,z_1^p-z_2^q)\in\ccc\{t,z_1,z_2\}^2$ be a~one-parameter deformation of $f(z_1,z_2)=(z_1^a+z_2^b,z_1^p-z_2^q)$. Assume that $a,b,p,q>1$ are integers such that $\mathrm{GCD}(p,q)=1$ and $bp<q$. Then $\m_0(F_{\underbar{t}})=bp$ for all $\underbar{t}\in\ccc$, i.e. $F$ is a~multiplicity-constant deformation. If $\underbar{t}\neq 0$ then $\ord F_{\underbar{t}}=1$ and we get $l_0(F_{\underbar{t}})=\m_0(F_{\underbar{t}})=bp$ by Lemma \ref{lem1.1}. Since $\ord F_0>1$ we get by the second part of Lemma \ref{lem1.1} that $l_0(F_0)<\m_0(F_0)=bp$.
\end{ex}
Note that C.~Bivi\`a-Ausina, see \cite{BA}, Corollary 2.5 proved a~result on the semicontinuity of the \L{}ojasiewicz exponent which however, does not imply our Theorem \ref{thm2.1}.

One can also indicate the deformations for which the \L{}ojasiewicz exponent is upper semicontinuous like multiplicity.
\begin{prop}\label{prop2.3}
Let $F\in\ccc\{t,z\}^n$ be a~deformation of $f\in\ccc\{z\}^n$ such that $\rank\left(\frac{\partial F_i}{\partial z_j}(\underbar{t},0)\right)\geqslant n-1$ for $\underbar{t}\in\ccc^k$ close to $0\in\ccc^k$. Then
$$
l_0(F_{\underbar{t}})\leqslant l_0(F_0)\quad\mbox{for }\underbar{t}\in\ccc^k\mbox{ close to }0.
$$
\end{prop}
{\sc Proof.} By Lemma \ref{lem1.1} we get $l_0(F_{\underbar{t}})=\m_0(F_{\underbar{t}})$ for $\underbar{t}\in\ccc^k$ close to $0$ and the proposition follows from the upper semicontinuity of the multiplicity.\hfill $\Box$
\begin{ex}\rm
Let $f(z)=(z_1^m,z_2,\ldots,z_n)$ with $m>1$ and let $F(t,z)=f(z_1+t,z_2,\ldots,z_n)-f(t,0,\ldots,0)=((z_1+t)^m-t^m,z_2,\dots,z_n)$ be a~one-parameter deformation of $f$. Then $F(t,z)$ satisfies the assumption of Proposition \ref{prop2.3}. Using Lemma \ref{lem1.1} we check that $l_0(F_{\underbar{t}})=\m_0(F_{\underbar{t}})=1$ for $\underbar{t}\neq 0$ and $l_0(F_0)=\m_0(F_0)=m$.
\end{ex}
In the example above the deformation of $f$ is given by the translation of coordinates. Even for such a~deformation the \L{}ojasiewicz exponent may be not upper semicontinuous:
\begin{ex}\label{ex2.5}\rm
Let $f(z_1,z_2,z_3)=(z_1^2,z_2^3,z_3^3-z_1z_2)\in\ccc\{z_1,z_2,z_3\}^3$ and let $F(t,z_1,z_2,z_3)=f(t+z_1,z_2,z_3)-f(t,0,0)=(2tz_1+z_1^2,z_2^3,-tz_2+z_3^3-z_1z_2)$. Then by Lemma \ref{lem1.1} we get $l_0(F_{\underbar{t}})=\m_0(F_{\underbar{t}})=9$ for $\underbar{t}\neq 0$. On the other hand $\m_0(F_0)=18$ and $l_0(F_0)=\frac{18}{5}$ (see Example \ref{ex3.5} of this note). The exponent $l_0(F_0)$ is attained on the arc $\phi(s)=(s^9,s^6,s^5)$.
\end{ex}
\begin{rem}\rm
The case of $\mu$-constant deformations of isolated hypersurface singularities is much more subtle. The Teissier's conjecture that ``$\mu$-constant implies the constancy of the \L{}ojasiewicz exponent'' \cite{T1}, Question on p. 278 is still open.
\end{rem}

\section{Characteristic polynomial and the \L{}o\-ja\-sie\-wicz exponent}\label{sect3}

Let $f=(f_1,\ldots,f_n)\in\ccc\{z\}^n$ be a~sequence of convergent power series defining a~finite holomorphic germ $f:(\ccc^n,0)\to(\ccc^n,0)$. Then the extension $\ccc\{z\}\supset\ccc\{f\}$ is a~finite $\ccc\{f\}$-module. For any $h\in\ccc\{z\}$ there is a~unique irreducible polynomial $Q_{f,h}=s^{m_h}+c_1(w)s^{m_h-1}+\cdots+c_{m_h}(w)\in\ccc\{w\}[s]$ in $n+1$ variables $(w,s)=(w_1,\ldots,w_n,s)$ such that $Q_{f,h}(f,h)=0$. It is called the minimal polynomial of $h$ relative to $f$. Its degree $m_{f,h}:=\deg_sQ_{f,h}$ divides the multiplicity $\m(f)$; we put $P_{f,h}=Q_{f,h}^r$, where $r=\frac{\m(f)}{m_{f,h}}$ and call $P_{f,h}$ the characteristic polynomial of $h$ relative to $f$. If $h(0)=0$ then $Q_{f,h}$ and consequently $P_{f,h}$ is a~distinguished polynomial.
\begin{rem}\rm Let $L=\ccc\{z\}_{(0)}$ and $K=\ccc\{f\}_{(0)}$ be fields of fractions of the ring $\ccc\{z\}$ and $\ccc\{f\}$, respectively. Then $Q_{f,h}(f,s)\in K[s]$ is the monic minimal polynomial of $h$ relative to the field extension $L/K$ and $P_{f,h}(f,s)$ is the characteristic polynomial of $h$ relative to $L/K$. For the various equivalent definitions of the characteristic polynomial, see Zariski-Samuel \cite{ZS}, Chapter II, \S{}10.
\end{rem} 
The lemma below follows immediately from the R\"uckert-Weierstrass parametrization theorem, see \cite{A}, \S{}31, (31.23).
\begin{lem}\label{lem3.2}
Let $P(w,s)=s^m+a_1(w)s^{m-1}+\cdots+a_m(w)\in\ccc\{w\}[s]$ be a~distinguished polynomial of degree $m=\m_0(f)$ and let $h\in\ccc\{z\}$, $h(0)=0$. Then the two conditions are equivalent
\begin{itemize}
\item[(i)] $P(w,s)$ is the characteristic polynomial of $h$ relative to $f$,
\item[(ii)] Let $U$ and $V$ be neighbourhoods of $0\in\ccc^n$ such that the mapping $U\ni\underbar{z}\to f(\underbar{z})\in V$ is a~$\m_0(f)$-sheeted branched covering and $h=h(\underbar{z})$ is convergent in $V$. Then the set $\{(\underbar{w},\underbar{s})\in V\times\ccc:P(\underbar{w},\underbar{s})=0\}$ is the image of $U$ by the mapping $U\ni\underbar{z}\to(f(\underbar{z}),h(\underbar{z}))\in V\times\ccc$, provided that $U$, $V$ are small enough.
\end{itemize}
\end{lem}

To study the \L{}ojasiewicz exponent $l_0(f)$ it is useful to consider the inequalities of the type\\[1.5ex]
(\mbox{\L{}})\hfill\qquad\qquad\qquad $|h(\underbar{z})|\leqslant c|f(\underbar{z})|^\theta\quad\mbox{near the origin }0\in\ccc^n.\hfill\mbox{}
$\\[1.5ex]
The least upper bound of the set of all $\theta>0$ for which (\L{}) holds for some constant $c>0$ in a~neighbourhood $U\subset\ccc^n$ of $0$ will be denoted $o_f(h)$ and called the \L{}ojasiewicz exponent of $h$ relative to $f$.
\begin{lem}\label{lem3.3}
Let $P_{f,h}(w,s)=s^m+a_1(w)s^{m-1}+\cdots+a_m(w)\in\ccc\{w,s\}$ be the characteristic polynomial of $h\in\ccc\{z\}$, $h\neq 0$, relative to $f$. Let $I=\{i\in\{1,\ldots,m\}:a_i\neq 0\}$. Then
$$
o_f(h)=\min_{i\in I}\left\{\frac{1}{i}\ord a_i\right\}.
$$
\end{lem}
{\sc Proof.} (after \cite{P2}, proof of Theorem 2.3). Let $U$ and $V$ be neighbourhoods of $0\in\ccc^n$ such that the mapping $U\ni\underbar{z}\to f(\underbar{z})\in V$ is an $\m_0(f)$-sheeted branched covering and $h=h(z)$ is convergent in $V$. Let $P(w,s)$ be the characteristic polynomial of $h$ relative to $f$. Then by Lemma \ref{lem3.2} we have that the inequality $|h(\underbar{z})|\leqslant c|f(\underbar{z})|^\theta$, $\underbar{z}\in U$, is equivalent to the estimate\\[1.5ex]
($*$)$\quad\qquad\{(\underbar{w},\underbar{s})\in V\times\ccc:P(\underbar{w},\underbar{s})=0\}\subset\{(\underbar{w},\underbar{s})\in V\times\ccc:|\underbar{s}|\leqslant|\underbar{w}|^\theta\}
$\\[1.5ex]
for $U$, $V$ small enough.

Let $\displaystyle\Theta_0=\min_{i\in I}\left\{\frac{1}{i}\ord a_i\right\}$. It is easy to check (see \cite{P0}, Proposition 2.2) that $\Theta_0$ is the largest number $\theta>0$ for which ($*$) holds. This proves the lemma.\hfill$\Box$
\begin{lem}\label{lem3.4}
$\displaystyle\qquad l_0(f)=\left(\min_{i=1}^n\{o_f(z_i)\}\right)^{-1}$.
\end{lem}
{\sc Proof.} Obvious.\hfill$\Box$

\begin{ex}\rm\label{ex3.5} Let us get back to Example \ref{ex2.5}. Let $f=(f_1,f_2,f_3)=(z_1^2,z_2^3,z_3^3-z_1z_2)$. We have $\m_0(f)=18$. The characteristic polynomials of $z_1$ and $z_2$ are $(s_1^2-w_1)^9$ and $(s_2^3-w_2)^6$ respectively, hence $o_f(z_1)=\frac{1}{2}$, $o_f(z_2)=\frac{1}{3}$. To calculate $o_f(z_3)$ let us observe that
$$
P(w,s)=(s^3-w_3)^6-w_1^3w_2^2
$$
is the characteristic polynomial of $h=z_3$ relative to $f$. Indeed, we have $P(f,z_3)=0$ in $\ccc\{z\}$ and $P(w,s)$ is irreducible: if $u$ is a~variable then $P(u,u,0,s)=s^{18}-u^5$ is irreducible, whence $P(w,s)$ is irreducible. 

Write $P(w,s)=s^{18}-6w_3s^{15}+\cdots+(w_3^6-w_1^3w_2^2)$. Using Lemma \ref{lem3.3} we check that $o_f(z_3)=\frac{\ord (w_3^6-w_1^3w_2^2)}{18}=\frac{5}{18}$. Then we get $l_0(f)=\left(\min\{\frac{1}{2},\frac{1}{3},\frac{5}{18}\}\right)^{-1}=\frac{18}{5}$.
\end{ex}
\begin{lem}\label{lem3.6}
Let $F=F(t,z)\in\ccc\{t,z\}^n$ be a~multiplicity-constant deformation of a~finite germ $f:(\ccc^n,0)\to(\ccc^n,0)$ and let $h\in\ccc\{z\}$. $h(0)=0$. Let $P_h(t,w,s)=s^m+a_1(t,w)s^{m-1}+\cdots+a_m(t,w)\in\ccc\{t,w\}[s]$ be the characteristic polynomial of $h$ relative to $(t,F(t,z))$. Then for $\underbar{t}\in\ccc^k$ close enough to $0\in\ccc^k$ the polynomial $P_h(\underbar{t},w,s)=s^m+a_1(\underbar{t},w)s^{m-1}+\cdots+a_m(\underbar{t},w)\in\ccc\{w\}[s]$ is the characteristic polynomial of $h$ relative to $F(\underbar{t},z)\in\ccc\{z\}^n$.
\end{lem}
{\sc Proof.} There exist arbitrary small neighbourhoods $U$ and $V$  of $0\in\ccc^n$ and $W$ of $0\in\ccc^k$ such that mapping $W\times U\ni(\underbar{t},\underbar{z})\to (\underbar{t},F(\underbar{t},\underbar{z}))\in W\times V$ is $\m_0(f)$-sheeted branched covering. Since $F=F(t,z)$ is a~multiplicity-constant deformation the mappings $U\ni\underbar{z}\to F(\underbar{t},\underbar{z})\in V$ for $\underbar{t}\in W$ are also $\m_0(f)$-sheeted branched coverings. Fix $h=h(z)\in\ccc\{z\}$, $h(0)=0$. Shrinking the neighbourhoods $W\times U$ and $W\times V$ we get by Lemma \ref{lem3.2} that the image of $W\times U$ under the mapping $W\times U\ni(\underbar{t},\underbar{z})\to (\underbar{t},F(\underbar{t},\underbar{z}),h(\underbar{z}))\in W\times V\times\ccc$ has the equation $P_h(t,w,s)=0$ in $W\times V\times\ccc$. Therefore the image of $U$ under the mapping $U\ni\underbar{z}\to (F(\underbar{t},\underbar{z}),h(\underbar{z}))\in V\times\ccc$ has the equation $P_h(\underbar{t},w,s)=0$ in $V\times\ccc$. Using again Lemma \ref{lem3.2} we have that $P_h(\underbar{t},w,s)$ is the characteristic polynomial of $h$ relative to $F(\underbar{t},z)$.\hfill$\Box$

\section{Proof of the main result}\label{sect4}

Let us begin with
\begin{thm}\label{thm4.1}
Let $F=F(t,z)\in\ccc\{t,z\}^n$ be a~multiplicity-constant deformation of a~finite germ $f:(\ccc^n,0)\to(\ccc^n,0)$. Let $h\in\ccc\{z\}$, $h\neq 0$. Then
$$
o_{F_{\underbar{t}}}(h)\leqslant o_{F_0}(h)\quad \mbox{for }\underbar{t}\in\ccc^k\mbox{ close to }0\in\ccc^k.
$$
Moreover, if $F$ is a~one-parameter deformation ($k=1$), then $o_{F_{\underbar{t}}}(h)$ is constant for $\underbar{t}\neq 0$ close to $0\in\ccc$.
\end{thm}
{\sc Proof.} Let $P_h(t,w,s)=s^m+a_1(t,w)s^{m-1}+\cdots+a_m(t,w)\in\ccc\{t,w\}[s]$ be the characteristic polynomial of $h$ relative to $(t,F(t,z))\in\ccc\{t,z\}^{k+n}$. Then by Lemma \ref{lem3.6} for $\underbar{t}\in\ccc^k$ close to $0\in\ccc^k$ we have that $P_h(\underbar{t},w,s)=s^m+a_1(\underbar{t},w)s^{m-1}+\cdots+a_m(\underbar{t},w)\in\ccc\{w\}[s]$ is the characteristic polynomial of $h$ relative to $F_{\underbar{t}}$. By Lemma \ref{lem3.3} $\displaystyle o_{F_{\underbar{t}}}(h)=\inf_i\left\{\frac{\ord a_i(\underbar{t},w)}{i}\right\}\leqslant\inf_i\left\{\frac{\ord a_i(0,w)}{i}\right\}=o_{F_0}(h)$ for $\underbar{t}\in\ccc^k$ close to $0\in\ccc^k$ since $\ord a_i(\underbar{t},w)\leqslant\ord a_i(0,w)$ if $|\underbar{t}|$ is small. If $k=1$ then $\ord a_i(\underbar{t},w)\equiv const$ for $\underbar{t}\neq 0$ close to $0\in\ccc$ and $o_{F_{\underbar{t}}}(h)=const$.\hfill$\Box$\vspace{1ex}

\noindent{\sc Proof of Theorem \ref{thm2.1}} Use Theorem \ref{thm4.1} and Lemma \ref{lem3.4}.\hfill$\Box$

\section{\L{}ojasiewicz exponent and the Newton polygon}

Let $P(w,s)=s^m+a_1(w)s^{m-1}+\cdots+a_m(w)\in\ccc\{w,s\}$ be a~distinguished polynomial in variables $(w,s)=(w_1,\ldots,w_n,s)$. Put $a_0(w)=1$ and $I=\{i:a_i\neq 0\}$. The Newton polygon $\mathcal{N}(P)$ of $P$ is defined to be
$$
\mathcal{N}(P)=\mathrm{convex}\bigcup_{i\in I}\left((\ord a_i,m-i)+\mathbb{R}_+^2\right),\quad\mbox{where }\mathbb{R}_+=\{a\in\mathbb{R}:a\geqslant 0\}.
$$
Then $\mathcal{N}(P)$ intersects the vertical axis at point $(0,m)$ and the horizontal axis at point $(\ord a_m,0)$ provided that $a_m\neq 0$. Note that $\displaystyle\theta(P):=\inf_i\left\{\frac{\ord a_i}{i}\right\}$ is equal to the inclination of the first side of the Newton polygon $\mathcal{N}(P)$, see \cite{T4}.

Let $f:(\ccc^n,0)\to(\ccc^n,0)$ be a~finite holomorphic germ and let $h\in\ccc\{z\}$, $h(0)=0$, $h\neq 0$ in $\ccc\{z\}$. We put
$$
\mathcal{N}(f,h)=\sigma(\mathcal{N}(P_{f,h})),
$$
where $\sigma$ is the symmetry of $\mathbb{R}_+^2$ given by $\sigma(\alpha,\beta)=(\beta, \alpha)$, and call $\mathcal{N}(f,h)$ the Newton polygon of $h$ relative to $f$.

>From the proof of Theorem \ref{thm4.1} it follows the semicontinuity of the Newton polygon in Teissier's sens, see \cite{T3}, pp. and \cite{T1}.
\begin{thm}
Let $F=F(t,z)\in\ccc\{t,z\}^n$ be a~multiplicity-constant deformation of $f$. Then
$$
\mathcal{N}(F_{\underbar{t}},h)\subset\mathcal{N}(F_0,h)\quad\mbox{for }\underbar{t}\in\ccc^k\mbox{ close to }0.
$$
If $k=1$ then $\mathcal{N}(F_{\underbar{t}}.h)$ does not depend on $\underbar{t}$ provided that $\underbar{t}\neq 0$ is close to $0\in\ccc$.
\end{thm}
Observe that $\mathcal{N}(f,h)$ intersects the horizontal axis at point $(\m_0(f),0)$. The intersection of the last edge (with vertex at $(\m_0(f),0)$) of $\mathcal{N}(f,h)$ is equal to $\frac{1}{o_f(h)}$. We will prove elsewhere that $\mathcal{N}(f,h)$ is identical to the Newton polygon of the pair of ideals $\I(f)$, $\I(h)=(h)\ccc\{z\}$ introduced by Teissier in \cite{T2}. In the notation of \cite{LJ-T}, Compl\'ement 2 we have $\mathcal{N}(f,h)=\mathcal{N}_{\I(f)}(h)$.

\noindent{}Kielce University of Technology\\
\noindent{}Department of Mathematics\\
Al. 1000 L PP 7\\
25-314 Kielce, Poland\\
E:mail: matap@tu.kielce.pl
\end{document}